\documentclass[12pt]{article}
\usepackage[latin1]{inputenc}
\usepackage[british]{babel}
\usepackage{cmap}
\usepackage{lmodern}

\usepackage{amssymb, amsmath, amsthm}
\usepackage[a4paper,top=25mm,bottom=25mm,left=25mm,right=25mm]{geometry}
\usepackage{ragged2e}

\usepackage{authblk} 
\usepackage{pifont}
\usepackage{graphicx}
\usepackage[usenames,dvipsnames,svgnames,table]{xcolor}
\usepackage[figuresright]{rotating}
\usepackage{xtab} 
\usepackage{longtable} 
\usepackage{multirow}
\usepackage{footnote}
\usepackage[stable]{footmisc}
\usepackage{chngpage} 
\usepackage{pdflscape} 
\usepackage[nottoc,notlot,notlof]{tocbibind} 

\usepackage{pgfplots}
\pgfplotsset{every tick label/.append style={font=\footnotesize}}
\pgfplotsset{compat=1.14}
\usepackage{setspace}

\usepackage{array}
\newcolumntype{K}[1]{>{\centering\arraybackslash$}p{#1}<{$}}

\makesavenoteenv{tabular}
\usepackage{tabularx}
\usepackage{booktabs}
\usepackage{threeparttable}
\usepackage[referable]{threeparttablex} 
\newcolumntype{R}{>{\raggedleft\arraybackslash}X}
\newcolumntype{L}{>{\raggedright\arraybackslash}X}
\newcolumntype{C}{>{\centering\arraybackslash}X}
\newcolumntype{A}{>{\columncolor{gray!25}}C}
\newcolumntype{a}{>{\columncolor{gray!25}}c}

\newlength{\tablen}

\usepackage{dcolumn} 
\newcolumntype{.}{D{.}{.}{-1}}

\usepackage{tikz}
\usetikzlibrary{arrows, calc, matrix, patterns, positioning, trees}
\usepackage[semicolon]{natbib}
\usepackage[hyphens]{url}
\usepackage{hyperref} 
\hypersetup{
  colorlinks   = true,    		
  urlcolor     = blue,    		
  linkcolor    = blue,    		
  citecolor    = ForestGreen	
}
\usepackage{microtype}
\usepackage[justification=centering]{caption} 

\usepackage[labelformat=simple]{subcaption}

\DeclareCaptionLabelFormat{parenthesis}{(#2)}
\makeatletter
\renewcommand\p@subfigure{\arabic{figure}.}
\makeatother

\DeclareCaptionLabelFormat{parenthesis}{(#2)}
\makeatletter
\renewcommand\p@subtable{\arabic{table}.}
\makeatother

\usepackage{enumitem}

\setlist[itemize]{leftmargin=2.5\parindent}
\setlist[enumerate]{leftmargin=2.5\parindent}

\def\addlegendimage{\csname pgfplots@addlegendimage\endcsname}

\theoremstyle{plain}

\theoremstyle{definition}

\newtheorem{definition}{Definition}[section]
\newtheorem{example}{Example}[section]

\theoremstyle{remark}

\newtheorem{remark}{Remark}

\makeatletter
\let\@fnsymbol\@alph
\makeatother

\def\keywords{\vspace{.5em} 
{\noindent \textit{Keywords}: }}

\def\AMS{\vspace{.5em} 
{\noindent \textbf{\emph{MSC} class}: }}

\def\JEL{\vspace{.5em} 
{\noindent \textbf{\emph{JEL} classification number}: }}

\title{Inconsistency thresholds for incomplete \\ pairwise comparison matrices}
\author{Kolos Csaba \'Agoston\thanks{~Email: \emph{kolos.agoston@uni-corvinus.hu} \newline Corvinus University of Budapest (BCE), Department of Operations Research and Actuarial Sciences, Budapest, Hungary}
$\qquad \qquad$
\href{https://sites.google.com/view/laszlocsato}{L\'aszl\'o Csat\'o}\thanks{~Corresponding author. Email: \emph{laszlo.csato@sztaki.hu} \newline
Institute for Computer Science and Control (SZTAKI), E\"otv\"os Lor\'and Research Network (ELKH), Laboratory on Engineering and Management Intelligence, Research Group of Operations Research and Decision Systems, Budapest, Hungary \newline
Corvinus University of Budapest (BCE), Department of Operations Research and Actuarial Sciences, Budapest, Hungary}} 
\date{\today}

\def\Dedication{
\begin{small}
{\noindent
``\emph{Mathematics is the part of physics where experiments are cheap.}''\footnote{~Source: \citet[p.~229]{Arnold1998}.}
}
\end{small}

\flushright
\begin{small}
(Vladimir Igorevich Arnold: \emph{On teaching mathematics})
\end{small}

\vspace{0.5cm} 
\justify }

\begin{document}

\maketitle
\thispagestyle{empty}
\Dedication

\begin{abstract}
\noindent
Pairwise comparison matrices are increasingly used in settings where some pairs are missing. However, there exist few inconsistency indices for similar incomplete data sets and no reasonable measure has an associated threshold. This paper generalises the famous rule of thumb for the acceptable level of inconsistency, proposed by Saaty, to incomplete pairwise comparison matrices. The extension is based on choosing the missing elements such that the maximal eigenvalue of the incomplete matrix is minimised. Consequently, the well-established values of the random index cannot be adopted: the inconsistency of random matrices is found to be the function of matrix size and the number of missing elements, with a nearly linear dependence in the case of the latter variable. Our results can be directly built into decision-making software and used by practitioners as a statistical criterion for accepting or rejecting an incomplete pairwise comparison matrix.



\AMS{90B50, 91B08}

\JEL{C44, D71}
\end{abstract}

\clearpage

\section{Introduction} \label{Sec1}

Pairwise comparisons form an essential part of many decision-making techniques, especially since the appearance of the popular \emph{Analytic Hierarchy Process (AHP)} methodology \citep{Saaty1977, Saaty1980}. Despite simplifying the issue to evaluating objects pair by pair, the tool of pairwise comparisons presents some challenges due to the possible lack of consistency: if alternative $A$ is two times better than alternative $B$ and alternative $B$ is three times better than alternative $C$, then alternative $A$ is not necessarily six times better than alternative $C$. The origin of similar inconsistencies resides in asking seemingly ``redundant'' questions. Nonetheless, additional information is often required to increase robustness \citep{SzadoczkiBozokiTekile2022}, and inconsistency usually does not cause a serious problem until it remains at a moderate level. 

Inconsistent preferences call for quantifying the level of inconsistency. The first and by far the most extensively used index has been proposed by the founder of the AHP, \emph{Thomas L.\ Saaty} \citep{Saaty1977}. He has also provided a sharp threshold to decide whether a pairwise comparison matrix has an acceptable level of inconsistency or not.

This widely accepted rule of inconsistency has been constructed for the case when all comparisons are known. However, there are at least three arguments why \emph{incomplete} pairwise comparisons should be considered in decision-making models \citep{Harker1987a}:
\begin{itemize}
\item
in the case of a large number $n$ of alternatives, completing all $n(n - 1)/2$ pairwise comparisons is resource-intensive and might require much effort from experts suffering from a lack of time;
\item
unwillingness to make a direct comparison between two alternatives for ethical, moral, or psychological reasons;
\item
the decision-makers may be unsure of some of the comparisons, for instance, due to limited knowledge on the particular issue.
\end{itemize}
In certain settings, both incompleteness and inconsistency are inherent features of the data. The beating relation in sports is rarely transitive and some players/teams have never played against each other \citep{BozokiCsatoTemesi2016, Csato2013a, Csato2017c, CsatoPetroczy2021b, ChaoKouLiPeng2018}. Analogously, there exists no guarantee for consistency when the pairwise comparisons are given by the bilateral remittances between countries \citep{Petroczy2021a}, or by the preferences of students between universities \citep{CsatoToth2020}.

Finally, note that pairwise comparison matrices are usually filled sequentially by the decision-makers, see e.g.\ the empirical research conducted by \citet{BozokiDezsoPoeszTemesi2013}. If the degree of inconsistency is monitored continuously during this process, the decision-maker might be warned immediately after the appearance of an unexpected value \citep{BozokiFulopKoczkodaj2011}. Consequently, there is a higher chance that the problem can be solved easily compared to the usual case when the supervision of the comparisons is asked only after all pairwise comparisons are given. This is especially important as these values are often provided by experts who suffer from a lack of time.

Let us see an example, where the missing elements are denoted by $\ast$:
\[
\mathbf{A} = \left[
\begin{array}{K{3em} K{3em} K{3em} K{3em}}
    1     		& 2  	& \ast    	& 4 \\
    1/2  	& 1       	& 2	& \ast \\
    \ast  		& 1/2 	& 1      	& 2 \\
    1/4  	&  \ast  	& 1/2 	& 1 \\
\end{array}
\right].
\]
Pairwise comparison matrix $\mathbf{A}$ is inconsistent because $a_{12} \times a_{23} \times a_{34} = 2 \times 2 \times 2 = 8 \neq 4 = a_{14}$. But it remains unknown whether this deviation can be tolerated or not.

The current paper aims to provide thresholds of acceptability for pairwise comparison matrices with missing entries. We want to follow the concept of Saaty as closely as possible. Therefore, the unknown elements are considered as variables to be chosen to reduce the inconsistency of the parametric complete pairwise comparison matrix, that is, to minimise its maximal eigenvalue as suggested by \citet{ShiraishiObataDaigo1998} and \citet{ShiraishiObata2002}. The main challenge resides in the calculation of the random index, a key component of Saaty's threshold: the optimal completion of each randomly generated incomplete pairwise comparison matrix should be found separately in order to obtain the minimal value of the Perron root of the completed matrix \citep{BozokiFulopRonyai2010}.

Inconsistency indices are thoroughly researched in the literature \citep{Brunelli2018}. There exist several attempts to calculate thresholds for Saaty's index under different assumptions \citep{AlonsoLamata2006, BozokiRapcsak2008, Ozdemir2005}, as well as for various inconsistency indices such as the geometric consistency index \citep{AguaronMoreno-Jimenez2003}, or the Salo--Hamalainen index \citep{AmentaLucadamoMarcarelli2020}. \citet{LiangBrunelliRezaei2019} propose consistency thresholds for the Best Worst Method (BWM).

On the other hand, the study of inconsistency indices for incomplete pairwise comparisons has been started only recently.
\citet{SzybowskiKulakowskiPrusak2020} introduce two new inconsistency measures based on spanning trees. \citet{KulakowskiTalaga2020} adapt several existing indices to analyse incomplete data sets but do not provide any threshold. To conclude, without the present contribution, one cannot decide whether the inconsistency of the above incomplete pairwise comparison matrix $\mathbf{A}$ is excessive or not.
Thus our work fills a substantial research gap.

Even though \citet{Forman1990} computes random indices for incomplete pairwise comparison matrices, his solution is based on the proposal of \citet{Harker1987a}. That introduces an auxiliary matrix for any incomplete pairwise comparison matrix instead of filling it by optimising an objective function as we do. Our approach is probably closer to Saaty's concept since the auxiliary matrix of \citet{Harker1987a} is not a pairwise comparison matrix.

The paper is structured as follows.
Section~\ref{Sec2} presents the fundamentals of pairwise comparison matrices and inconsistency measures. Incomplete pairwise comparison matrices and the eigenvalue minimisation problem are introduced in Section~\ref{Sec3}. Section~\ref{Sec4} discusses the details of computing the random index. The inconsistency thresholds are reported in Section~\ref{Sec5}. A numerical example is provided in Section~\ref{Sec6}, and a real life application in Section~\ref{Sec7}. Finally, Section~\ref{Sec8} offers a summary and directions for future research.

\section{Pairwise comparison matrices and inconsistency} \label{Sec2}

The pairwise comparisons of the alternatives are collected into a matrix $\mathbf{A} = \left[ a_{ij} \right]$ such that the entry $a_{ij}$ is the numerical answer to the question ``How many times alternative $i$ is better than alternative $j$?''
Let $\mathbb{R}_+$ denote the set of positive numbers, $\mathbb{R}^n_+$ denote the set of positive vectors of size $n$ and $\mathbb{R}^{n \times n}_+$ denote the set of positive square matrices of size $n$ with all elements greater than zero, respectively.

\begin{definition} \label{Def21}
\emph{Pairwise comparison matrix}:
Matrix $\mathbf{A} = \left[ a_{ij} \right] \in \mathbb{R}^{n \times n}_+$ is a \emph{pairwise comparison matrix} if $a_{ji} = 1 / a_{ij}$ for all $1 \leq i,j \leq n$.
\end{definition}

Let $\mathcal{A}$ denote the set of pairwise comparison matrices and $\mathcal{A}^{n \times n}$ denote the set of pairwise comparison matrices of size $n$, respectively.

\begin{definition} \label{Def22}
\emph{Consistency}:
A pairwise comparison matrix $\mathbf{A} = \left[ a_{ij} \right] \in \mathcal{A}^{n \times n}$ is \emph{consistent} if $a_{ik} = a_{ij} a_{jk}$ for all $1 \leq i,j,k \leq n$. Otherwise, it is said to be \emph{inconsistent}.
\end{definition}


According to the famous Perron--Frobenius theorem, for any pairwise comparison matrix $\mathbf{A} \in \mathcal{A}$, there exists a unique positive weight vector $\mathbf{w}$ satisfying $\mathbf{A} \mathbf{w} = \lambda_{\max}(\mathbf{A}) \mathbf{w}$ and $\sum_{i=1}^n w_i = 1$, where $\lambda_{\max}(\mathbf{A})$ is the maximal or Perron eigenvalue of matrix $\mathbf{A}$.

Saaty has considered an affine transformation of this eigenvalue.

\begin{definition} \label{Def24}
\emph{Consistency index}:
Let $\mathbf{A} \in \mathcal{A}^{n \times n}$ be any pairwise comparison matrix of size $n$.
Its \emph{consistency index} is
\[
CI(\mathbf{A}) = \frac{\lambda_{\max}(\mathbf{A}) - n}{n-1}.
\]
\end{definition}

Since $CI(\mathbf{A}) = 0 \iff \lambda_{\max}(\mathbf{A}) = n$, the consistency index $CI$ is a reasonable measure of how far a pairwise comparison matrix is from a consistent one \citep{Saaty1977, Saaty1980}.
\citet{AupetitGenest1993} provide a tight upper bound for the value of $CI$ when the entries of the pairwise comparison matrix are expressed on a bounded scale.

Saaty has recommended using a discrete scale for the matrix elements, i.e., for all $1 \leq i,j \leq n$:
\begin{equation} \label{Saaty_scale}
a_{ij} \in \left\{ 1/9,\, 1/8,\, 1/7, \dots ,\, 1/2,\, 1,\, 2, \dots ,\, 8,\, 9 \right\}.
\end{equation}

A normalised measure of inconsistency can be obtained as suggested by Saaty.

\begin{definition} \label{Def25}
\emph{Random index}:
Consider the set $\mathcal{A}^{n \times n}$ of pairwise comparison matrices of size $n$.
The corresponding \emph{random index} $RI$ is provided by the following algorithm \citep{AlonsoLamata2006}:
\begin{itemize}
\item
Generating a large number of pairwise comparison matrices such that each entry above the diagonal is drawn independently and uniformly from the Saaty scale \eqref{Saaty_scale}.
\item
Calculating the consistency index $CI$ for each random pairwise comparison matrix.
\item
Computing the mean of these values.
\end{itemize}
\end{definition}

\begin{table}[t!]
\caption{The values of the random index for complete pairwise comparison matrices}
\centering
\label{Table1}
    \begin{tabularx}{\textwidth}{l CCCCCCC} \toprule
    Matrix size & 4     & 5     & 6     & 7     & 8     & 9     & 10 \\ \midrule
    Random index $RI_n$ & 0.884 & 1.109 & 1.249 & 1.341 & 1.404 & 1.451 & 1.486 \\ \bottomrule
    \end{tabularx}
\end{table}

Several authors have published slightly different random indices depending on the simulation method and the number of generated matrices involved, see \citet[Table~1]{AlonsoLamata2006}.
The random indices $RI_n$ are reported in Table~\ref{Table1} for $4 \leq n \leq 10$ as provided by \citet{BozokiRapcsak2008} and validated by \citet{CsatoPetroczy2021a}. These estimates are close to the ones given in previous works \citep{AlonsoLamata2006, Ozdemir2005}. \citet[Table~3]{BozokiRapcsak2008} uncovers how $RI_n$ depends on the largest element of the ratio scale.

\begin{definition} \label{Def26}
\emph{Consistency ratio}:
Let $\mathbf{A} \in \mathcal{A}^{n \times n}$ be any pairwise comparison matrix of size $n$. Its \emph{consistency ratio} is $CR(\mathbf{A}) = CI(\mathbf{A}) / RI_n$.
\end{definition}

Saaty has proposed a threshold for the acceptability of inconsistency, too.

\begin{definition} \label{Def27}
\emph{Acceptable level of inconsistency}:
Let $\mathbf{A} \in \mathcal{A}^{n \times n}$ be any pairwise comparison matrix of size $n$.
It is sufficiently close to a consistent matrix and therefore can be accepted if $CR(\mathbf{A}) \leq 0.1$.
\end{definition}

Even though applying a crisp decision rule on the fuzzy concept of ''large inconsistency'' is strange \citep{Brunelli2018} and there exist sophisticated statistical studies to test consistency \citep{LinKouErgu2013, LinKouErgu2014}, it is assumed throughout the paper that the 10\% rule is a well-established standard worth generalising to incomplete pairwise comparison matrices.

\section{The eigenvalue minimisation problem for incomplete pairwise comparison matrices} \label{Sec3}

Certain entries of a pairwise comparison matrix are sometimes missing.

\begin{definition} \label{Def31}
\emph{Incomplete pairwise comparison matrix}:
Matrix $\mathbf{A} = \left[ a_{ij} \right]$ is an \emph{incomplete pairwise comparison matrix} if $a_{ij} \in \mathbb{R}_+ \cup \{ \ast \}$ such that for all $1 \leq i,j \leq n$, $a_{ij} \in \mathbb{R}_+$ implies $a_{ji} = 1 / a_{ij}$ and $a_{ij} = \ast$ implies $a_{ji} = \ast$.
\end{definition}

Let $\mathcal{A}_{\ast}^{n \times n}$ denote the set of incomplete pairwise comparison matrices of size $n$.

The graph representation of incomplete pairwise comparison matrices is a convenient tool to visualise the structure of known elements.

\begin{definition} \label{Def32}
\emph{Graph representation}:
An incomplete pairwise comparison matrix $\mathbf{A} \in \mathcal{A}_{\ast}^{n \times n}$ can be represented by the undirected graph $G = (V,E)$, where the vertices $V = \left\{ 1,2, \dots ,n \right\}$ correspond to the alternatives and the edges in $E$ are associated with the known matrix entries outside the diagonal, that is, $e_{ij} \in E \iff a_{ij} \neq \ast \text{ and } i \neq j$.
\end{definition}

To summarise, there are no edges for the missing elements ($a_{ij} = \ast$) as well as for the entries of the diagonal ($a_{ii}$).

In the case of an incomplete pairwise comparison matrix $\mathbf{A}$, \citet{ShiraishiObataDaigo1998} and \citet{ShiraishiObata2002} consider an eigenvalue optimisation problem by substituting the $m$ missing elements of matrix $\mathbf{A}$ above the diagonal with positive values arranged in the vector $\mathbf{x} \in \mathbb{R}^m_+$, while the reciprocity condition is preserved:
\begin{equation} \label{EV_optimisation}
\min_{\mathbf{x} \in \mathbb{R}^m_+} \lambda_{\max} \left( \mathbf{A}(\mathbf{x}) \right).
\end{equation}
The motivation is clear, all missing entries should be chosen to get a matrix that is as close to a consistent one as possible in terms of the consistency index $CI$.

According to \citet[Section~3]{BozokiFulopRonyai2010}, \eqref{EV_optimisation} can be transformed into a convex optimisation problem. The authors also give the necessary and sufficient condition for the uniqueness of the solution: the graph $G$ representing the incomplete pairwise comparison matrix $\mathbf{A}$ should be connected. This is an intuitive and almost obvious requirement since the relation of two alternatives cannot be established if they are not compared at least indirectly, through other alternatives.

\section{The calculation of the random index for incomplete pairwise comparison matrices} \label{Sec4}

Consider an incomplete pairwise comparison matrix $\mathbf{A} \in \mathcal{A}_{\ast}^{n \times n}$ and a complete pairwise comparison matrix $\mathbf{B} \in \mathcal{A}^{n \times n}$, where $b_{ij} = a_{ij}$ if $a_{ij} \neq \ast$.
Let $\mathbf{A}(\mathbf{x}) \in \mathcal{A}^{n \times n}$ be the optimal completion of $\mathbf{A}$ according to \eqref{EV_optimisation}. Clearly, $\lambda_{\max} \left( \mathbf{A}(\mathbf{x}) \right) \leq \lambda_{\max}(\mathbf{B})$, hence $CI \left( \mathbf{A}(\mathbf{x}) \right) \leq CI(\mathbf{B})$.
This implies that the value of the random index $RI_n$, calculated for complete pairwise comparison matrices, cannot be applied in the case of an incomplete pairwise comparison matrix because its consistency index $CI$ is obtained through optimising (i.e.\ minimising) its level of inconsistency.

Consequently, by adopting the numbers from Table~\ref{Table1}, the ratio of incomplete pairwise comparison matrices with an acceptable level of inconsistency will exceed the concept of Saaty and this discrepancy increases as the number of missing elements grows. In the extreme case when graph $G$ is a spanning tree of a complete graph with $n$ nodes (thus it is a connected graph consisting of exactly $n-1$ edges without cycles), the corresponding incomplete matrix can be filled out such that consistency is achieved.

Therefore, the random index needs to be recomputed for incomplete pairwise comparison matrices, and its value will supposedly be a monotonically decreasing function of $m$, the number of missing elements.

\begin{remark} \label{Rem1}
In the view of the Saaty scale \eqref{Saaty_scale}, there are at least three different ways to choose the missing entries $x_k$, $1 \leq k \leq m$:
\begin{enumerate}
\item \label{Method1}
\emph{Method 1}: $x_k \in \mathbb{R}_+$, namely, each missing entry can be an arbitrary positive number;
\item \label{Method2}
\emph{Method 2}: $1/9 \leq x_k \leq 9$, namely, the missing entries cannot be higher (lower) than the theoretical maximum (minimum) of the known elements;
\item \label{Method3}
\emph{Method 3}: $x_k \in \left\{ 1/9,\, 1/8,\, 1/7, \dots ,\, 1/2,\, 1,\, 2, \dots ,\, 8,\, 9 \right\}$, namely, each missing entry is drawn from the discrete Saaty scale.
\end{enumerate}
\end{remark}

Let us illustrate the three approaches listed in Remark~\ref{Rem1}.

\begin{example} \label{Examp1}
Take the following incomplete pairwise comparison matrix:
\[
\mathbf{A} = \left[
\begin{array}{K{3em} K{3em} K{3em} K{3em}}
    1     & \ast  & 9     & \ast \\
    \ast  & 1     & 2     & 8 \\
     1/9  &  1/2  & 1     & 4 \\
    \ast  &  1/8  &  1/4  & 1 \\
\end{array}
\right].
\]

\begin{figure}[t!]
\centering
\begin{tikzpicture}[scale=1, auto=center, transform shape, >=triangle 45]
\tikzstyle{every node}=[draw,shape=circle];
  \node (n1) at (135:2) {$1$};
  \node (n2) at (45:2)  {$2$};
  \node (n3) at (315:2) {$3$};
  \node (n4) at (225:2) {$4$};

  \foreach \from/\to in {n1/n3,n2/n3,n2/n4,n3/n4}
    \draw (\from) -- (\to);
\end{tikzpicture}

\caption{The graph representation of the pairwise comparison matrix $\mathbf{A}$ in Example~\ref{Examp1}}
\label{Fig1}
\end{figure}
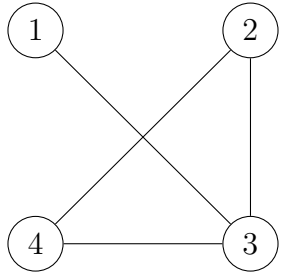

The corresponding undirected graph $G$ is depicted in Figure~\ref{Fig1}. Note that $G$ would be a spanning tree without the edge between nodes $2$ and $4$ and $a_{24} = 8 = 2 \times 4 = a_{23} a_{34}$. Consequently, $\mathbf{A}$ can be filled out consistently in a unique way:
\[
\mathbf{A^1} = \left[
\begin{array}{K{3em} K{3em} K{3em} K{3em}}
    1     & \mathbf{9/2}  & 9     & \mathbf{36} \\
    \mathbf{2/9}  & 1     & 2     & 8 \\
     1/9  &  1/2  & 1     & 4 \\
    \mathbf{1/36}  &  1/8  &  1/4  & 1 \\
\end{array}
\right].
\]

The first technique (Method~\ref{Method1} in Remark~\ref{Rem1}) results in $\mathbf{A^1}$ with $\lambda_{\max} \left( \mathbf{A^1} \right) = 4$.

On the other hand, $\mathbf{A^1}$ is not valid under Method~\ref{Method2} in Remark~\ref{Rem1} because $a_{14}^1 = 36 > 9$, that is, the consistent filling is not allowed as being outside the Saaty scale \eqref{Saaty_scale}. The optimal complete pairwise comparison matrix $\mathbf{A^2}$ is given by the solution of the convex eigenvalue minimisation problem \eqref{EV_optimisation} with the additional constraints $1/9 \leq x_k \leq 9$ for all $1 \leq k \leq m$ and is as follows:
\[
\mathbf{A^2} = \left[
\begin{array}{K{3em} K{3em} K{3em} K{3em}}
    1    &  \mathbf{9/4}  & 9     & \mathbf{9} \\
    \mathbf{4/9}  & 1     & 2     & 8 \\
    1/9  &  1/2  & 1     & 4 \\
    \mathbf{1/9}  &  1/8  &  1/4  & 1 \\
\end{array}
\right],
\]
where $\lambda_{\max} \left( \mathbf{A^2} \right) = 4.1855$.

Finally, $\mathbf{A^2}$ is not valid under Method~\ref{Method3} in Remark~\ref{Rem1} because $a_{12}^2 = 9/4 \notin \mathbb{Z}$, that is, even though the optimal filling by Method~\ref{Method2} does not contain any value exceeding the bounds of the Saaty scale \eqref{Saaty_scale}, some of them are not integers or the reciprocals of integers. Hence, the best possible filling on the Saaty scale \eqref{Saaty_scale} is
\[
\mathbf{A^3} = \left[
\begin{array}{K{3em} K{3em} K{3em} K{3em}}
    1     & \mathbf{2}  & 9     & \mathbf{9} \\
    \mathbf{1/2}  & 1     & 2     & 8 \\
     1/9  &  1/2  & 1     & 4 \\
    \mathbf{1/9}  &  1/8  &  1/4  & 1 \\
\end{array}
\right],
\]
which leads to $\lambda_{\max} \left( \mathbf{A^3} \right) = 4.1874$.
\end{example}

Among the three ideas in Remark~\ref{Rem1}, Method~\ref{Method1} always leads to the smallest dominant eigenvalue, followed by Method~\ref{Method2}, whereas Method~\ref{Method3} provides the greatest optimum of problem \eqref{EV_optimisation} as can be seen from the restrictions in Remark~\ref{Rem1}.

We implement Method~\ref{Method2} to calculate the random indices $RI_n$.
The first reason is that the algorithm for the $\lambda_{\max}$-optimal completion \citep[Section~5]{BozokiFulopRonyai2010} involves an exogenously given tolerance level to determine how accurate are the coordinates of the eigenvector associated with the dominant eigenvalue as a stopping criterion. Consequently, it cannot be chosen appropriately if the matrix entries and the elements of the weight vector can differ substantially: the consistent completion of an incomplete pairwise comparison matrix with $n$ alternatives may contain $(1/9)^{(n-1)}$ or $9^{(n-1)}$ as an element if the corresponding graph is a chain. Furthermore, it remains questionable why elements below or above the Saaty scale \eqref{Saaty_scale} are allowed for the missing entries if they are prohibited in the case of known elements.
On the other hand, Method~\ref{Method3} presents a discrete optimisation problem that is more difficult to handle than its continuous analogue of Method~\ref{Method2}.
To summarise, since the process is based on generating a large number of random incomplete pairwise comparison matrices to be filled out optimally, it is necessary to reduce the complexity of optimisation problem \eqref{EV_optimisation} by using Method~\ref{Method2}.

A complete pairwise comparison matrix of size $n$ can be represented by a complete graph where the degree of each node is $n-1$. Hence, the graph corresponding to an incomplete pairwise comparison matrix is certainly connected if $m \leq n-2$, implying that the solution of the $\lambda_{\max}$-optimal completion is unique.
However, the graph might be disconnected if $m \geq n-1$, in which case it makes no sense to calculate the consistency index $CI$ of the incomplete pairwise comparison matrix.
Furthermore, if $m > n(n-1)/2 - (n-1)$, then there are less than $n-1$ known elements, and the graph is always disconnected.

If the number of missing entries is exactly $m = n(n-1)/2 - (n-1) = (n-1)(n-2)/2$, then the graph is connected if and only if it is a spanning tree. Even though these incomplete pairwise comparison matrices certainly have a consistent completion under Method~\ref{Method1}, this does not necessarily hold under Method~\ref{Method2} when the missing entries cannot be arbitrarily large/small.

\section{Generalised thresholds for the consistency ratio} \label{Sec5}

As we have argued in Section~\ref{Sec4}, the value of the random index $RI_{n,m}$ probably depends not only on the size $n$ of the incomplete pairwise comparison matrix but on the number of its missing elements $m$, too.
Thus the random index is computed according to the following procedure (cf.\ Definition~\ref{Def25}): 
\begin{enumerate}
\item \label{Step1}
Generating an incomplete pairwise comparison matrix $\mathbf{A}$ of size $n$ with $m$ missing entries above the diagonal such that each element above the diagonal is drawn independently and uniformly from the Saaty scale \eqref{Saaty_scale}, while the place of the unknown elements above the diagonal is chosen randomly.
\item
Checking whether the graph $G$ representing the incomplete pairwise comparison matrix $\mathbf{A}$ is connected or disconnected.
\item
If graph $G$ is connected, optimisation problem \eqref{EV_optimisation} is solved by the algorithm for the $\lambda_{\max}$-optimal completion \citep[Section~5]{BozokiFulopRonyai2010} with restricting all entries in $\mathbf{x} \in \mathbb{R}^m_+$ according to Method~\ref{Method2} in Remark~\ref{Rem1} to obtain the minimum of $\lambda_{\max} \left( \mathbf{A(x)} \right)$ and the corresponding complete pairwise comparison matrix $\mathbf{\hat{A}}$.
\item \label{Step4}
Computing and saving the consistency index $CI \left( \mathbf{\hat{A}} \right)$ based on Definition~\ref{Def24}.
\item
Repeating Steps~\ref{Step1}--\ref{Step4} to get 1 million random matrices with a connected graph representation, and calculating the mean of the consistency indices $CI$ from Step~\ref{Step4}.
\end{enumerate}

\begin{table}[t!]
\caption{The values of the random index \\ for incomplete pairwise comparison matrices}
\centering
\label{Table2}
\begin{threeparttable}
\rowcolors{3}{gray!20}{}
    \begin{tabularx}{\textwidth}{c CCCC} \toprule \hiderowcolors    
    \multirow{2}{*}{Missing elements $m$} & \multicolumn{4}{c}{Matrix size $n$} \\
          & 4     & 5     & 6     & 7     \\ \bottomrule \showrowcolors
    0     & 0.884 & 1.109 & 1.249 & 1.341 \\
    1     & 0.583 (0.531) & 0.925 (0.485) & 1.128 (0.400) & 1.256 (0.330) \\
    2     & 0.306 (0.387) & 0.739 (0.452) & 1.007 (0.392) & --- \\
    3     & 0.053 (0.073) & 0.557 (0.405) & 0.883 (0.380) & --- \\
    4     & ---   & 0.379 (0.340) & 0.758 (0.364) & --- \\
    5     & ---   & 0.212 (0.247) & 0.634 (0.344) & --- \\
    6     & ---   & 0.059 (0.068) & 0.510 (0.317) & --- \\
    7     & ---   & ---   & 0.389 (0.281) & --- \\
    8     & ---   & ---   & 0.271 (0.234) & --- \\
    9     & ---   & ---   & 0.161 (0.170) & --- \\ \toprule
    \end{tabularx}
\begin{tablenotes} \footnotesize
\item
All values are based on 1 million matrices. Standard deviations are given in parenthesis.
\end{tablenotes}
\end{threeparttable}
\end{table}

Our central result is reported in Table~\ref{Table2}, which is an extension of Table~\ref{Table1} to the case when some pairwise comparisons are unknown. The values in the first row, which coincide with the ones from Table~\ref{Table1}, confirm the integrity of the proposed technique to compute the thresholds for the consistency index $CI$.
The role of missing elements cannot be ignored at all commonly used significance levels as reinforced by the t-test: for any given $n$, the values of $RI_{n,m}$ are statistically different from each other.

\begin{figure}[t!]

\begin{tikzpicture}
\begin{axis}[
name = axis1,
title = {Matrix size $n=5$},
title style = {font=\small},
xlabel = Number of missing elements $m$,
x label style = {font=\small},
width = 0.5\textwidth,
height = 0.4\textwidth,
ymajorgrids = true,
xmin = 0,
xmax = 6,
] 
\addplot [red, thick, mark=pentagon*] coordinates {
(0,1.109)
(1,0.9246087)
(2,0.7387896)
(3,0.556883)
(4,0.37861)
(5,0.211943)
(6,0.0591422)
};
\end{axis}

\begin{axis}[
at = {(axis1.south east)},
xshift = 0.1\textwidth,
title = {Matrix size $n=6$},
title style = {font=\small},
xlabel = Number of missing elements $m$,
x label style = {font=\small},
width = 0.5\textwidth,
height = 0.4\textwidth,
ymajorgrids = true,
xmin = 0,
xmax = 9,
]
\addplot [blue, thick, dashdotted, mark=asterisk, mark options={solid,thick}] coordinates {
(0,1.249)
(1,1.127994)
(2,1.007017)
(3,0.8827498)
(4,0.7582752)
(5,0.6343154)
(6,0.50963)
(7,0.388586)
(8,0.2713317)
(9,0.160653)
};
\end{axis}
\end{tikzpicture}

\caption{The random index $RI_{n,m}$ as the function of the number of missing entries $m$}
\label{Fig2}
\end{figure}
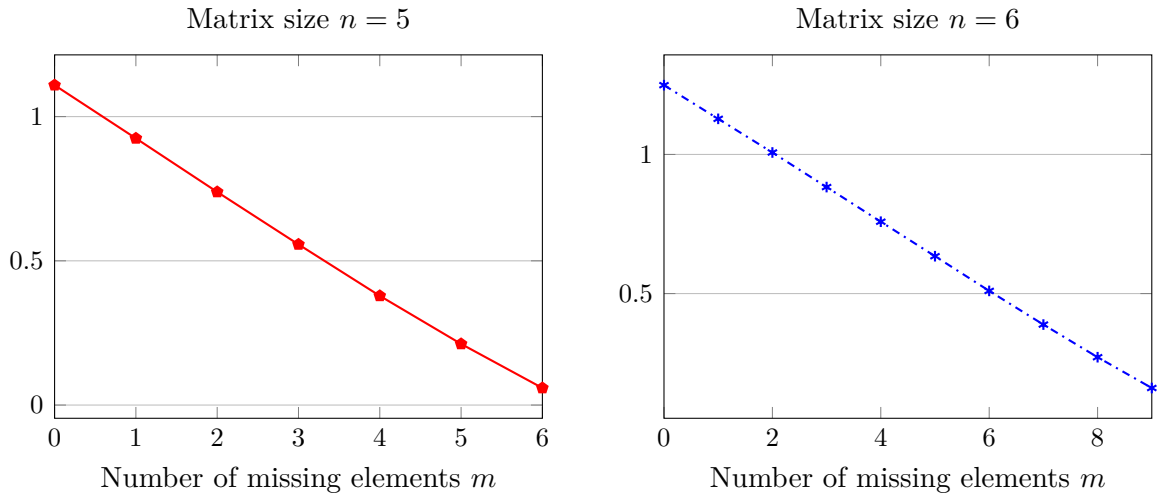


Recall that the maximal number of missing elements is at most $n(n-1)/2 - (n-1) = (n-1)(n-2)/2$ if connectedness is not violated, and this value is $3$ if $n=4$, 6 if $n=5$, and 10 if $n=6$. Some thresholds are lacking from Table~\ref{Table2}---for example, the pair $n=7$ and $m=4$---due to excessive computation time ($>48$ hours).

However, $RI_{n,m}$ can be easily predicted as follows.
Figure~\ref{Fig2} reveals that the random index is monotonically decreasing as the function of missing values $m$ according to common intuition.
Furthermore, the dependence is nearly linear, thus a plausible estimation is provided by the below formula, which requires only the ``omnipresent'' Table~\ref{Table1}:
\begin{equation} \label{RI_approx}
RI_{n,m} \approx \left[ 1 - \frac{2m}{(n-1)(n-2)} \right] RI_{n,0}.
\end{equation}
Obviously, \eqref{RI_approx} returns $RI_{n,0}$ if there are no missing elements ($m=0$). On the other hand, $m = (n-1)(n-2)/2$ means that the graph representing the incomplete pairwise comparison matrix is either unconnected, or it is a spanning tree, thus the matrix can be filled consistently if there is no restriction on its elements. Formula \eqref{RI_approx} immediately follows by assuming a linear function for intermediate values of $m$.

\begin{table}[t!]
\caption{Approximation of the random index for incomplete \\ pairwise comparison matrices according to equation \eqref{RI_approx}}
\centering
\label{Table3}
\rowcolors{3}{gray!20}{}
    \begin{tabularx}{\textwidth}{ccCc} \toprule \hiderowcolors
    \multirow{2}{*}{Matrix size $n$} & \multirow{2}{*}{Missing elements $m$} & \multicolumn{2}{c}{Value of $RI_{n,m}$} \\
          &       & Computed & Approximated by formula \eqref{RI_approx} \\ \bottomrule \showrowcolors
    7     & 4     & 0.998 & 0.983 \\
    8     & 5     & 1.088 & 1.070 \\
    9     & 6     & 1.158 & 1.140 \\
    10    & 7     & 1.215 & 1.197 \\ 
    15    & 4     & 1.519 & 1.514 \\
    15    & 8     & 1.453 & 1.445 \\
    15    & 12    & 1.387 & 1.375 \\ \toprule
    \end{tabularx}
\end{table}

According to the ``case studies'' in Table~\ref{Table3}, \eqref{RI_approx} gives at least a reasonable guess of $RI_{n,m}$ without much effort, even though it somewhat underestimates the true value. The discrepancy is mainly caused by $RI_{n,(n-1)(n-2)/2}$ being larger than zero (see Table~\ref{Table2}) as incomplete pairwise comparison matrices represented by a spanning tree can be made consistent only if the missing elements can be arbitrary, but not if they are bounded to the interval $\left[ 1/9,\, 9 \right]$.

Definition~\ref{Def26} can be modified straightforwardly to derive the consistency ratio $CR$ for any incomplete pairwise comparison matrix.

\begin{definition} \label{Def51}
\emph{Consistency ratio}:
Let $\mathbf{A} \in \mathcal{A}_{\ast}^{n \times n}$ be any incomplete pairwise comparison matrix of size $n$ with $m$ missing entries above the diagonal and $\mathbf{\hat{A}}$ be the complete pairwise comparison matrix given by the optimal filling of $\mathbf{A}$.
The \emph{consistency ratio} of the incomplete matrix $\mathbf{A}$ is $CR(\mathbf{A}) = CI(\mathbf{\hat{A}}) / RI_{n,m}$.
\end{definition}

The popular 10\% threshold of Definition~\ref{Def27} can be adopted without any changes.

In the applications of the AHP methodology, the optimal number of alternatives does not exceed nine \citep{SaatyOzdemir2003}. Random indices for complete pairwise comparison matrices have been determined for $n \leq 16$ in \citet{AguaronMoreno-Jimenez2003} and for $n \leq 15$ in \citet{AlonsoLamata2006}. The corresponding thresholds for incomplete pairwise comparison matrices can be calculated offline by a supercomputer and built into any software used by practitioners. If these are not available, formula \eqref{RI_approx} provides a good approximation for any number of alternatives $n$ and missing elements $m$, see Table~\ref{Table3}.

\section{An illustrative example} \label{Sec6}

In this section, we highlight the implications of the calculated thresholds for the random index by a numerical illustration. It has been chosen to be simple but expressive. With three alternatives and one missing entry, the matrix can be filled out consistently. Therefore, the number of alternatives is four. Again, there exists a consistent filling if there are three missing elements, hence their number is two. Furthermore, they are in different rows, which is the more likely case.

\begin{example} \label{Examp2}
Take the following parametric incomplete pairwise comparison matrix of size $n=4$ with $m=2$ missing elements:
\[
\mathbf{A}(\alpha, \beta) = \left[
\begin{array}{K{3em} K{3em} K{3em} K{3em}}
    1     		& \alpha  	& \ast    	& \beta \\
    1/\alpha  	& 1       	& \alpha	& \ast \\
    \ast  		& 1/\alpha 	& 1      	& \alpha \\
    1/\beta  	&  \ast  	& 1/\alpha 	& 1 \\
\end{array}
\right].
\]

\begin{table}[t!]
\caption{Consistency indices of the parametric incomplete \\ pairwise comparison matrix $\mathbf{A}(\alpha, \beta)$ in Example~\ref{Examp2}}
\centering
\label{Table4}
\rowcolors{3}{gray!20}{}
\centerline{
\begin{threeparttable}
    \begin{tabularx}{1.125\textwidth}{c CCCCC CCCC} \toprule \hiderowcolors
    \multirow{2}{*}{Value of $\beta$} & \multicolumn{9}{c}{Value of $\alpha$} \\
          &  1/5  &  1/4  &  1/3  &  1/2  & 1     & 2     & 3     & 4     & 5 \\ \bottomrule \showrowcolors
     1/9  & 0.1495 & \textit{0.0818} & \textbf{0.0253} & \textbf{0.0003} & 0.1031 & 0.4187 & 0.7338 & 1.0344 & 1.3214 \\
     1/8  & 0.1637 & 0.0921 & \textbf{0.0311} & \textbf{0} & 0.0921 & 0.3940 & 0.6982 & 0.9890 & 1.2671 \\
     1/7  & 0.1807 & 0.1047 & \textit{0.0383} & \textbf{0.0004} & \textit{0.0805} & 0.3670 & 0.6592 & 0.9393 & 1.2076 \\
     1/6  & 0.2015 & 0.1202 & \textit{0.0477} & \textbf{0.0017} & \textit{0.0680} & 0.3374 & 0.6160 & 0.8842 & 1.1414 \\
     1/5  & 0.2278 & 0.1401 & \textit{0.0601} & \textbf{0.0046} & \textit{0.0547} & 0.3042 & 0.5673 & 0.8220 & 1.0667 \\
     1/4  & 0.2624 & 0.1667 & \textit{0.0774} & \textbf{0.0100} & \textit{0.0404} & 0.2663 & 0.5113 & 0.7500 & 0.9801 \\
     1/3  & 0.3114 & 0.2048 & 0.1031 & \textbf{0.0201} & \textbf{0.0253} & 0.2217 & 0.4444 & 0.6637 & 0.8759 \\
     1/2  & 0.3891 & 0.2663 & 0.1462 & \textit{0.0404} & \textbf{0.0100} & 0.1667 & 0.3599 & 0.5536 & 0.7426 \\
    1     & 0.5476 & 0.3940 & 0.2394 & 0.0921 & \textbf{0} & 0.0921 & 0.2394 & 0.3940 & 0.5476 \\
    2     & 0.7426 & 0.5536 & 0.3599 & 0.1667 & \textbf{0.0100} & \textit{0.0404} & 0.1462 & 0.2663 & 0.3891 \\
    3     & 0.8759 & 0.6637 & 0.4444 & 0.2217 & \textbf{0.0253} & \textbf{0.0201} & 0.1031 & 0.2048 & 0.3114 \\
    4     & 0.9801 & 0.7500 & 0.5113 & 0.2663 & \textit{0.0404} & \textbf{0.0100} & \textit{0.0774} & 0.1667 & 0.2624 \\
    5     & 1.0667 & 0.8220 & 0.5673 & 0.3042 & \textit{0.0547} & \textbf{0.0046} & \textit{0.0601} & 0.1401 & 0.2278 \\
    6     & 1.1414 & 0.8842 & 0.6160 & 0.3374 & \textit{0.0680} & \textbf{0.0017} & \textit{0.0477} & 0.1202 & 0.2015 \\
    7     & 1.2076 & 0.9393 & 0.6592 & 0.3670 & \textit{0.0805} & \textbf{0.0004} & \textit{0.0383} & 0.1047 & 0.1807 \\
    8     & 1.2671 & 0.9890 & 0.6982 & 0.3940 & 0.0921 & \textbf{0} & \textbf{0.0311} & 0.0921 & 0.1637 \\
    9     & 1.3214 & 1.0344 & 0.7338 & 0.4187 & 0.1031 & \textbf{0.0003} & \textbf{0.0253} & \textit{0.0818} & 0.1495 \\ \toprule
    \end{tabularx}
\begin{tablenotes} \footnotesize
\item
\textbf{Bold} numbers indicate that the consistency ratio $CR \left( \mathbf{\hat{A}}(\alpha, \beta) \right) = CI \left( \mathbf{\hat{A}} (\alpha, \beta) \right) / RI_{4,2}$ is below the 10\% threshold.
\item
\textit{Italic} numbers indicate that $CI \left( \mathbf{\hat{A}} (\alpha, \beta) \right) / RI_{4,0}$ is below the 10\% threshold but the consistency ratio $CR \left( \mathbf{\hat{A}} (\alpha, \beta) \right) = CI \left( \mathbf{\hat{A}} (\alpha, \beta) \right) / RI_{4,2}$ is above it.
\end{tablenotes}
\end{threeparttable}
}
\end{table}

Now $RI_{4,0} \approx 0.884$ and $RI_{4,2} \approx 0.356$ from Table~\ref{Table2}. There are three instances where the optimal filling of matrix $\mathbf{A}(\alpha, \beta)$ results in a consistent pairwise comparison matrix:
\[
(\alpha,\, \beta) \in \left\{ \left( \frac{1}{2},\, \frac{1}{8} \right); \left( 1,\, 1 \right); \left( 2,\, 8 \right) \right\}.
\]
They should be accepted under any circumstances.

Examine what happens if $\alpha = 1$ is fixed. Then $\beta = 3$ implies $CI \left( \mathbf{\hat{A}}(1,3) \right) \approx 0.0253 < 0.1 \times RI_{4,2}$, which still corresponds to an acceptable level of inconsistency. However, $CI \left( \mathbf{\hat{A}}(1,4) \right) \approx 0.0404 > 0.1 \times RI_{4,2}$, making it necessary to reduce inconsistency if $\beta = 4$. On the other hand, $CI \left( \mathbf{\hat{A}}(1,4) \right) \approx 0.0404 < 0.1 \times RI_{4,0}$, thus the optimally filled out incomplete pairwise comparison matrix might be accepted according to the ``standard'' threshold for complete matrices because the latter does not take into account the automatic reduction of inconsistency due to the optimisation procedure.

Table~\ref{Table4} reports the consistency index $CI$ of matrix $\mathbf{A}(\alpha, \beta)$ for some parameters $\alpha$ and $\beta$. $\alpha$ is restricted between $1/5$ and $5$ because $a_{12}(\alpha, \beta) \times a_{23}(\alpha, \beta) \times a_{34}(\alpha, \beta) = \alpha^3$ but $a_{14}(\alpha, \beta) = \beta$. Bold numbers correspond to the cases when inconsistency can be tolerated based on the approximated thresholds of Table~\ref{Table2}, while italic numbers show instances that can be accepted only if the optimal solution $\mathbf{A(x)}$ of \eqref{EV_optimisation} is considered as a (complete) pairwise comparison matrix and the threshold of 10\% is used for $CI \left( \mathbf{A(x)} \right) / RI_{4,0}$.
\end{example}

\section{A real life application: continuous monitoring of inconsistency} \label{Sec7}

\citet{BozokiDezsoPoeszTemesi2013} carried out a controlled experiment, where university students were divided into subgroups to make pairwise comparisons from different types of problems, with different number of alternatives in different questioning orders. Consequently, not only the complete pairwise comparison matrices are known but their incomplete submatrices obtained after a given number of comparisons was asked. We have picked one interesting matrix from this dataset.

\begin{example} \label{Examp3}
The following pairwise comparison matrix reflects the opinion of a decision-maker on how much more a summer house is liked compared to another summer house on a numerical scale:
\[
\mathbf{A} = \left[
\begin{array}{K{3em} K{3em} K{3em} K{3em} K{3em} K{3em}}
    1     & 2     & 7     & 7     & 7     & 5 \\
     1/2  & 1     & 5     & 7     & 5     & 2 \\
     1/7  &  1/5  & 1     & 1     &  1/5  &  1/5 \\
     1/7  &  1/7  & 1     & 1     &  1/3  &  1/3 \\
     1/7  &  1/5  & 5     & 3     & 1     & 3 \\
     1/5  &  1/2  & 5     & 3     &  1/3  & 1 \\
\end{array}
\right].
\]
The questioning order of the $15$ comparisons is $a_{12}$, $a_{64}$, $a_{51}$, $a_{32}$, $a_{56}$, $a_{13}$, $a_{24}$, $a_{61}$, $a_{43}$, $a_{52}$, $a_{14}$, $a_{35}$, $a_{26}$, $a_{45}$, and $a_{36}$. This procedure, proposed by \citet{Ross1934}, optimises two objective functions: it maximises the distances for the same alternatives to reappear and aims to balance the number of the first and second positions in the comparison for every alternative.

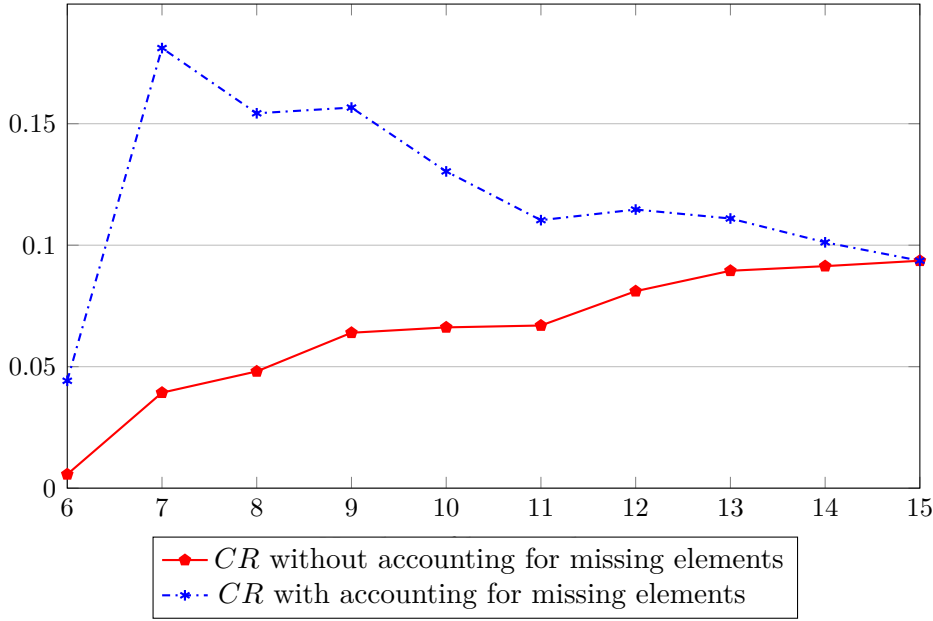
\begin{figure}[t!]
\centering
\begin{tikzpicture}
\begin{axis}[
name = axis1,
title style = {font=\small},
xlabel = Number of known elements,
x label style = {font=\small},
width = 0.8\textwidth,
height = 0.5\textwidth,
ymajorgrids = true,
xmin = 6,
xmax = 15,
ymin = 0,
yticklabel style = {scaled ticks=false,/pgf/number format/fixed,/pgf/number format/precision=2}, 
legend style = {font=\small,at={(0.1,-0.1)},anchor=north west,legend columns=1},
legend entries = {$CR$ without accounting for missing elements,$CR$ with accounting for missing elements$\quad \,$},
] 
\addplot [red, thick, mark=pentagon*] coordinates {
(6,0.00569675820656525)
(7,0.0392998983186549)
(8,0.048054261008807)
(9,0.0639635228182546)
(10,0.0661744251401121)
(11,0.0669433706965572)
(12,0.0810705892714171)
(13,0.0894788678943154)
(14,0.0913774395516413)
(15,0.093606)
};
\addplot [blue, thick, dashdotted, mark=asterisk, mark options={solid,thick}] coordinates {
(6,0.0441941)
(7,0.1811276)
(8,0.1542925)
(9,0.1566479)
(10,0.1303657)
(11,0.1103064)
(12,0.114674)
(13,0.1109822)
(14,0.1011795)
(15,0.093606)
};
\end{axis}
\end{tikzpicture}

\caption{The consistency ratio $CR$ as the function of known elements in Example~\ref{Examp3}}
\label{Fig3}
\end{figure}


Figure~\ref{Fig3} shows how inconsistency changes as more and more comparisons are given by the decision-maker. Following \citet[Figure~2]{BozokiDezsoPoeszTemesi2013}, the solid red line uses the random index associated with a complete $6 \times 6$ pairwise comparison matrix, which is not a valid approach according to Section~\ref{Sec4}. On the other hand, the dashed blue line is obtained by the values of the random index according to our computations, see Table~\ref{Table2}. The na\"ive approach indicates no problem around inconsistency, its level remains below the 10\% threshold during the filling in process. However, accounting for the number of missing elements reveals that inconsistency is substantially increased when the seventh comparison ($a_{24}$) is made. Even though the complete pairwise comparison matrix can be accepted with respect to inconsistency, continuous monitoring warns the decision-maker that this particular comparison is worth reconsidering.
\end{example}

\section{Conclusions} \label{Sec8}

The paper reports approximated thresholds for the most popular measure of inconsistency, proposed by Saaty, in the case of incomplete pairwise comparison matrices. They are determined by the value of the random index, that is, the average consistency ratio of a large number of random pairwise comparison matrices with missing elements. The calculation is far from trivial since a separate convex optimisation problem should be solved for each matrix to find the optimal filling of unknown entries. Numerical results uncover that the threshold depends not only on the size of the pairwise comparison matrix but on the number of missing entries, too. A plausible linear estimation of the random index has also been provided.

According to Table~\ref{Table2} and two examples, the extended values of the random index become indispensable in order to generalise Saaty's concept to incomplete comparisons. The associated thresholds can be directly programmed into decision-making software.

With the suggested rule of acceptability, the decision-maker can decide for any incomplete pairwise comparison matrix whether there is a need to revise earlier assessments or not. It allows the level of inconsistency to be monitored even before all comparisons are given, which may immediately indicate possible mistakes and suspicious entries. Therefore, the preference revision process can be launched as early as possible. It will be examined in future studies how this opportunity can be built into the known inconsistency reduction methods \citep{AbelMikhailovKeane2018, BozokiFulopPoesz2015, ErguKouPengShi2011, XuXu2020}.

\section*{Acknowledgements}
\addcontentsline{toc}{section}{Acknowledgements}
\noindent
We are grateful to \emph{S\'andor Boz\'oki}, \emph{Lavoslav {\v C}aklovic}, and \emph{Zsombor Sz\'adoczki} for useful advice. \\
Four anonymous reviewers provided valuable comments and suggestions on earlier drafts. \\
The research was supported by the MTA Premium Postdoctoral Research Program grant PPD2019-9/2019.

\bibliographystyle{apalike}
\bibliography{All_references}

\end{document}